\documentclass[12pt]{amsart} 
\usepackage{graphics}
\newfont{\sheaf}{eusm10 scaled\magstep1}

\def\Bbb{\bf}

\newcommand{\C}{\ensuremath{\mathbb{C}}}
\newcommand{\R}{\ensuremath{\mathbb{R}}} 
\newcommand{\Z}{\ensuremath{\mathbb{Z}}}

\newcommand{\D}{\ensuremath{\mathcal{D}}}
\newcommand{\MM}{\ensuremath{\mathcal{M}}}

\newcommand{\hol}{\ensuremath{\mathcal{O}}}
 
\newcommand{\PP}{\ensuremath{\mathbb{P}}} 

\newcommand{\T}{\ensuremath{\mathcal{T}}}
\newcommand{\BB}{\ensuremath{\mathcal{B}}}

\newcommand{\ra}{\ensuremath{\rightarrow}} 
\newcommand{\SSS}{\ensuremath{\mathcal{S}}} 
\newcommand{\PPP}{\ensuremath{\mathcal{P}}} 
 
\newtheorem{teo}{Theorem}[section] 
\newtheorem{df}[teo]{Definition}

\newtheorem{oss}[teo]{Remark} 

\newtheorem{ex}[teo]{Example}

\def\eea{\end{eqnarray*}}
\def\bea{\begin{eqnarray*}}

\def\eea{\end{eqnarray*}}
\def\bea{\begin{eqnarray*}}
\begin{document}

\title{Trecce, Mapping class group, fibrazioni di Lefschetz ed applicazioni al
diffeomorfismo di superficie algebriche.}

\author{Fabrizio Catanese}

\date{April 6, 2004}
\maketitle

\section{Introduzione}

Uno dei problemi fondamentali nella teoria delle superficie algebriche
e' quello di capire gli spazi di moduli delle superficie di tipo
generale, in particolare le loro componenti connesse, che
parametrizzano classi di equivalenza per deformazione di superficie
minimali di tipo generale.

Per un classico teorema di Ehresmann, due varieta' complesse equivalenti per
deformazione sono tra loro diffeomorfe.

Negli anni 80,   Friedman e Morgan (\cite{f-m1}) hanno congetturato che 
due superficie algebriche siano equivalenti per deformazione se e 
solamente se sono tra loro diffeomorfe tramite un diffeomorfismo che 
conserva l' orientazione (congettura da noi abbreviata con l' acronimo Def= Diff).

Dopo i primi contresempi di M.  Manetti (cf. \cite{man4}), dei contresempi piu' semplici
sono stati trovati dall' autore e da  Kharlamov-Kulikov (\cite{cat4} e \cite{k-k}). 

In entrambi questi ultimi lavori, vengono utilizzate coppie di superficie 
complesse coniugate.
Si puo' pero' correggere la congettura, e richiedere da un lato
un diffeomorfismo $ \phi : S' \ra S$ che porta la prima classe di Chern 
 $c_1(K_S) \in H^2 (S, \Z)$ del divisore canonico in $c_1(K_{S'})$ 
 (a seguito di risultati di Witten e Morgan si sa oggi che ogni
 diffeomorfismo porta  $c_1(K_S)$ o in  $c_1(K_{S'})$ o in  $- c_1(K_{S'})$).

D' altro lato si puo' restringere la congettura al caso di superficie 
semplicemente connesse, per cui il tipo topologico, per il famoso teorema di 
M. Freedman (\cite{free}), e'
univocamente determinato dai due invarianti numerici $K^2_S$ e $\chi (S):= \chi (\hol_S)$
e dalla classe di resto modulo  2 del numero $r (S) := $ divisibilita' di 
$c_1(K_S) $ in $H^2 (S, \Z)$.

In effetti, gli esempi di Manetti  hanno un primo gruppo di omologia non triviale
e gli altri esempi citati possiedono dei rivestimenti non ramificati 
di ordine finito per cui vale   Def = Diff .
 
Scopo fondamentale della conferenza e' stato quello di introdurre certi
concetti e tecniche, in parte classici, che sono stati utilizzati da  Bronislaw Wajnryb 
e dall' autore per dimostrare esplicitamente (\cite{c-w}) il diffeomorfismo di alcune
superficie algebriche assai elementari, ma che portano ai desiderati contresempi
alla congettura Def= Diff nel caso semplicemente connesso.

Alcune parole chiave sono : somma connessa,  
 chirurgie, gruppo  fondamentale di Poincare', gruppo delle trecce, superficie di Riemann
(curve complesse), gruppo del Dehn delle classi di applicazioni (Mapping class group).

Un aspetto  interessante ed ancora da esplorare e' la  sinergia fra la geometria
complessa e la geometria simplettica.

Questa sinergia ha come strumento chiave, a seguito di risultati fondamentali 
dovuti a Kas, Donaldson e Gompf (\cite{kas}, \cite{don6},\cite{gompf1}),
la estensione al caso di 4-varieta' simplettiche della teoria delle fibrazioni di Lefschetz,
che ha gia' svolto un ruolo centrale nello studio della topologia delle
varieta' algebriche.

 Mentre la dimostrazione del risultato nel lavoro dell' autore con Wajnryb e' assai
elaborata (e la sua parte piu' topologica  si articola anche attraverso 
una lunga serie di figure atte a mostrare la isotopia di certe curve), il suo spirito
e' assai generale ed abbiamo la speranza da una parte che il metodo
si presti a dimostrare il diffeomorfismo di importanti classi di superficie algebriche,
dall' altra che una variante del metodo porti anche a stabilire 
la equivalenza simplettica di certe superficie algebriche munite della struttura
simplettica canonica definita dall' autore in \cite{cat6}.

\tableofcontents

\section{Le idee di Riemann: varieta' differenziabili e 
topologiche ed operazioni su di esse.}

Ricordiamo la definizione di varieta', concetto introdotto essenzialmente da
B. Riemann nella sua dissertazione inaugurale per la  Libera Docenza.

Una varieta' $M$ di dimensione reale $m$ 
e' uno spazio topologico tale che per ogni punto $ p \in M$ 
esiste un aperto $ U_p$ contenente $p$ ed un omeomorfismo (detto carta locale)
$ \psi_p : U_p \ra V_p \subset \R^m$ su di un aperto $V_p$ di $\R^m$
ed inoltre vale la proprieta' fondamentale che,
nel dominio di definizione,

$\psi_{p'} \circ \psi_p^{-1}$ sia un:

\begin{itemize}
\item
Omeomorfismo (sull' immagine) per   $M$ {\bf varieta' topologica }
\item
Diffeomorfismo (sull' immagine) per   $M$ {\bf varieta' differenziabile }
\item
Biolomorfismo (sull' immagine) per   $M$ {\bf varieta' complessa } 
(in questo ultimo caso $m = 2n$, $\R^m = \C^n$). 

\end{itemize}

\begin{df}
La operazione di somma connessa $ M_1 \sharp M_2$ si puo' effettuare
nell' ambito delle varieta' differenziabili o topologiche della stessa dimensione.

Si scelgono rispettivi punti $p_i \in M_i$ e carte locali estendibili in un intorno
di $\overline{U_{p_i} }$ 
$$ \psi_{p_i}: U_{p_i}  \ra \cong B(0, \epsilon_i):
= \{ x \in \R^m | |x| < \epsilon_i \},$$
e, ponendo $B(0, \epsilon_i)^* := B(0, \epsilon_i) - \{ 0\} $,
 si considera un  diffeomorfismo 
$\psi : B(0, \epsilon_1)^* \ra  B(0, \epsilon_2)^*$ tale che 
$\psi (x) =  f (|x|) \frac{x}{|x|} $ per una opportuna funzione differenziabile
decrescente tale che $ f(0) = \epsilon_1 , f (\epsilon_2) = 0$.

Allora la somma connessa $ M_1 \sharp M_2$ e' lo spazio quoziente della unione 
disgiunta $ (M_1 - \{p_1\} ) \cup^o  (M_2 - \{p_2 \})$ per la relazione di 
equivalenza che identifica 
$ y \in (U_{p_1} - \{p_1\} ) $ a $  w \in (U_{p_2} - \{p_2\} ) $ se e solo se
$$ w = \psi_{p_2}^{-1} \circ \psi \circ \psi_{p_1} (y).$$ 
\end{df}

Sussiste il seguente

\begin{teo}
Il risultato della operazione di somma connessa e' ben definito, 
cioe' e' indipendente dalle scelte fatte dei due punti $p_1, p_2$ e 
della funzione scalare $f$. 
\end{teo}

\begin{ex}
L' esempio piu' intuitivo (vedi Figura 1) e' quello di due superficie 
di Riemann compatte (orientabili) $ M_1, M_2$ di rispettivi generi $g_1, g_2$:
allora $ M_1 \sharp M_2$ e' compatta orientabile ed ha genere $g_1 + g_2$. 
\end{ex}

\begin{figure}[htbp]
\begin{center}
\scalebox{1}{\includegraphics{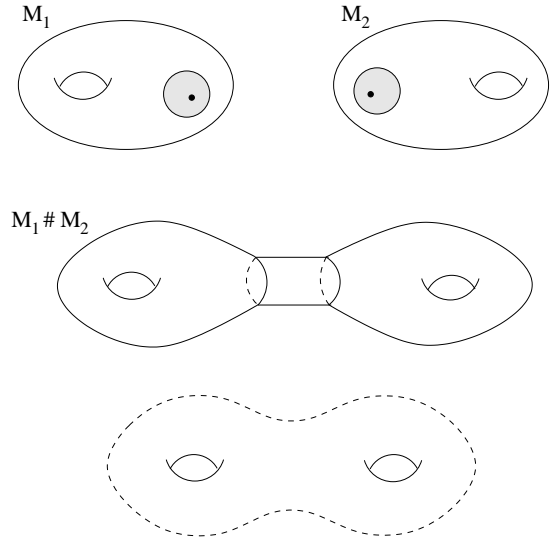}}
\end{center}
\caption{La Somma Connessa}
\label{figura1}
\end{figure}

\begin{oss}
L' operazione di somma connessa fra varieta' della stessa dimensione ci da' un
semigruppo: vale cioe' la associativita', ed abbiamo come elemento neutro la
sfera $ S^m : = \{ x \in \R^{m+1}| |x| = 1\}$. E' facile infatti verificare che
$ M \sharp S^m \cong M$. 
\end{oss}

\begin{df}
La varieta' $M$ si dice {\bf irriducibile} se $ M \cong M_1 \sharp  M_2$ implica che
$M_1$ oppure $M_2$ e' omotopicamente equivalente ad una sfera $S^m$. 
\end{df}

In dimensione due, ogni varieta' compatta di genere $ g \geq 1$ e' somma connessa di
$g$ tori (superficie di Riemann di genere $1$, diffeomorfi a $\R^2 / \Z^2$), e quindi
e' riducibile se $ g \geq 2$. Ma in dimensione piu' alta vi sono molte varieta' irriducibili,
ed allora si ricorre ad altri procedimenti costruttivi. Ecco un esempio ulteriore

\begin{df}
{\bf (CHIRURGIA)}

 Sia, per $i=1,2$,   $ N_i \subset  M_i$ una sottovarieta' differenziabile.
Questo vuol dire che per ogni $p \in M_i$ esiste una carta dove l' immagine 
di $N_i$ e' l' intersezione con un sottospazio di $\R^m$ dove si annullano certe 
coordinate, ed implica, per via del teorema delle funzioni implicite, che esiste
un aperto $U_i \supset N_i$ che e' diffeomorfo al fibrato normale $\nu_{N_i}$ della immersione
di $N_i \ra M_i$, tramite un diffeomorfismo che porta $N_i$ nella sezione nulla
di $\nu_{N_i}$. 

Supponiamo ora che esista un diffeomorfismo $\phi : N_1 \ra N_2$, ed un diffeomorfismo
$\psi : (\nu_{N_1} - N_1)\ra (\nu_{N_2} - N_2)$ compatibile con le proiezioni $p_i : \nu_{N_i} \ra N_i$
(cioe' tale che $p_2 \circ \psi = p_1$), e con la proprieta' di rovesciare la orientazione
sulle fibre. Si puo' allora prendere la varieta' $ M_1 \sharp_{\psi} M_2$, definita di nuovo
come quoziente della unione 
disgiunta $ (M_1 - N_1 ) \cup^o  (M_2 - N_2)$ per la relazione di 
equivalenza che identifica $(U_1 - N_1 )$ a $(U_2 - N_2 )$ tramite il diffeomorfismo
indotto da $\psi$.
\end{df}

\begin{oss}
Questa volta il risultato della operazione dipende dalla scelta di 
$\phi$ e di $\psi$.
\end{oss}

Il problema che nasce, dopo che abbiamo introdotto cosi' tanti modi di
 produrre nuove varieta', e' quello di riconoscere quando due varieta' sono tra 
loro omeomorfe o diffeomorfe.

Il metodo fondamentale della topologia algebrica e' quello di associare ad ogni
spazio topologico $X$ un oggetto algebrico $H(X)$ in modo funtoriale, cioe' in modo tale
che ad ogni applicazione continua $ f: X \ra Y$ corrisponda un unico omomorfismo 
$ H(f) : H(X) \ra H(Y)$. Un esempio archetipico e' dato dal gruppo fondamentale, introdotto
da Henri Poincare'.

\section{Il gruppo di Poincare' $\pi_1(X)$ .}

Sia $X$ uno spazio topologico e $p \in X$: allora il Gruppo Fondamentale (di Poincare')
di $X$ e' definito tramite l' insieme delle classi di omotopia
$$ \pi_1 (X, p) : = \{ \sigma : [0,1] \ra X | \sigma e' continua,  \sigma(0)= \sigma(1)=p 
\} / OMOTOPIA.$$
Dove, dati  spazi $Y, X$ e sottospazi $Y' \subset Y, X' \subset X$ 
(qui $Y = [0,1], Y' = \{ 0,1\}, X' = \{p\}$), due applicazioni continue 
$f,g: (Y, Y' ) \ra (X,X')$ (cioe' tali che portano $Y'$ in $X'$) si dicono
 OMOTOPE con estremi in $X'$ se  esiste una applicazione continua
$F :  (Y \times [0,1] , Y' \times [0,1]) \ra (X,X')$ tale che
$f(y) = F(y,0), g(y) = F(y,1)$.

La struttura di gruppo e' fornita dalla composizione di due cammini (a velocita' doppia),
che  si definisce  cosi' :
\begin{itemize}
\item
$ \sigma \circ \tau (t) : = \sigma(2t)$ per $ 0 \leq t \leq 1/2$, mentre
\item
 $ \sigma \circ \tau (t) : = \tau(2t-1)$ per $ 1/2 \leq t \leq 1$. 
\end{itemize}

Mentre il cammino inverso e' cosi' definito

$\sigma^{-1} (t) : = \sigma(t-1)$.

Nel caso di una varieta' connessa $M$ il gruppo fondamentale e' ben definito a prescindere
dalla scelta del punto base $p$, e quindi si denotera' sovente solo
con $\pi_1 (M)$.

\begin{ex}
Consideriamo la superficie di Riemann di genere $0$, ovvero la varieta' complessa
data dalla retta proiettiva $\PP^1_{\C} \cong S^2$, e prendiamo per $M$ il complementare
di un insieme di $k$ punti distinti: $ M : = \PP^1_{\C} - \{ p_1, \dots p_k\}$.

Poiche' a meno di diffeomorfismo possiamo sempre supporre che i punti dati
siano i punti $ p_h := - h$, e possiamo prendere come punto base il punto $p:=$
$ -  2 k i = - 2 k \surd \overline{-1}$, si definiscono i cammini $\gamma_h$ come la
composizione $ L_h \circ \beta_h L_h^{-1}$, $L_h$ essendo il segmento tra
$ -  2 k i$ e $p_h := -h$ percorso fino al punto $q_h$ su questo segmento che dista 
$1/10$ da $p_h$, e $\beta_h$ essendo la circonferenza di centro $p_h$ e raggio $1/10$, 
percorsa una volta in senso antiorario con punto di partenza in $q_h$.

Con queste notazioni abbiamo che il gruppo fondamentale e' il gruppo con generatori
e relazioni :
$$< \gamma_1, \dots \gamma_k | \Pi_{h=1}^k \gamma_h = 1 >. $$
I generatori descritti sono un esempio di una {\bf quasi-base geometrica}, 
che si ottiene
piu' generalmente scegliendo una serie di cammini tra il punto base $p$
ed i punti $p_h$,
che non si intersecano ne' si  auto-intersecano, e si susseguono 
l'uno dopo l'altro in 
ordine antiorario ( vedi Figura 2). 
\end{ex}

\begin{figure}[htbp]
\begin{center}
\scalebox{1}{\includegraphics{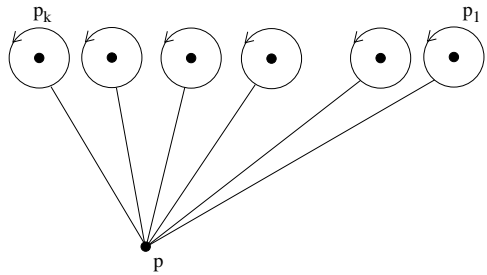}}
\end{center}
\caption{Una quasi-base geometrica}
\label{figura2}
\end{figure}

\begin{ex}
Nel caso di una varieta' complessa $C$ di dimensione $1$, 
o piu' generalmente di una varieta' topologica di dimensione 
due, il gruppo fondamentale
e' libero (caso particolare: l' esempio precedente, dove il gruppo e' liberamente
generato da $ \gamma_1, \dots \gamma_{k-1} $, che non soddisfano 
ad alcuna relazione)  se $C$ non e' compatta.

Se invece $C$ e' compatta  ed orientabile, allora $C$ e' ottenuta come somma
connessa di $g$ tori, ed il suo gruppo fondamentale ha la presentazione
$$< \alpha_1, \beta_1,  \dots   \alpha_g, \beta_g | \prod_{i=1}^g 
\alpha_i \beta_i \alpha_i^{-1} \beta_i^{-1} = 1 >.$$
Tale presentazione corrisponde alla ben nota realizzazione di $C$ come
spazio quoziente di un poligono con $2g$ lati, dove i lati sono identificati
secondo lo schema che alla successione di  lati percorsi in senso antiorario 
fa corrispondere la successione di cammini
 $\alpha_1, \beta_1, \alpha_1^{-1}, \beta_1^{-1}, \dots 
\alpha_g, \beta_g, \alpha_g^{-1}, \beta_g^{-1}$.
\end{ex}

\begin{oss}
Se prendiamo una somma connessa di due varieta' topologiche (che, si 
presuppone sempre tacitamente, siano connesse) $M_1$ ed $M_2$ di dimensione 
$m \geq 3$, allora il gruppo fondamentale della somma connessa
$M_1 \sharp M_2$ e' il prodotto libero 
$\pi_1(M_1) * \pi_1(M_2)$ di $\pi_1(M_1)$ e di $\pi_1(M_2)$.

In altre parole, date presentazioni
 $\pi_1(M_1) = <a_1, \dots a_h | R_1(a_i)=1, \dots R_t(a_i)=1> $ e
$\pi_1(M_2) = <a'_1, \dots a'_k | R'_1(a'_i)=1, \dots R'_s(a'_i)=1>$,
allora 
$$\pi_1(M_1) * \pi_1(M_2) =$$ 
$$= <a_1, \dots a_h,a'_1, \dots a'_k  | 
R_1(a_i)=1, \dots R_t(a_i)=1,R'_1(a'_i)=1, \dots R'_s(a'_i)=1>. $$

Questa osservazione fa vedere che l'operazione di somma connessa
fra varieta' non puo' essere una operazione di gruppo, poiche'
se $A$ e' un gruppo non banale, e $B$ e' un qualsiasi gruppo,
allora il gruppo $A$ e' contenuto nel prodotto libero $A*B$ e quindi
$A*B$ e' anche esso non banale, a differenza del gruppo fondamentale
dell' elemento neutro $S^m$, che \'e banale per $ m \geq 2$.
\end{oss}
Essendo il gruppo fondamentale un gruppo non Abeliano, esso da' una ampia
possibilita' per distinguere varieta' fra loro non omeomorfe.
Si sa, come osservato in \cite{seif-threl} che ogni gruppo finitamente
presentato e' il gruppo fondamentale di una varieta' topologica compatta di
dimensione $m \geq 4$, mentre assai difficile e' rispondere alla domanda:
quali sono i gruppi fondamentali di superficie complesse compatte (superficie =
varieta' di dimensione complessa 2)?

Come vedremo in seguito, esistono molti tipi topologici e differenziabili distinti
di varieta' semplicemente connesse, cioe' con un gruppo fondamentale banale.

\section{ Trecce e Mapping class group.}

Si deve ad Emil Artin (cf. \cite{art1}, ed altri articoli cf. \cite{art})
la elegante definizione del gruppo delle trecce, che
fornisce uno strumento potente, ma assai complicato da maneggiare, per lo
studio della topologia differenziale delle varieta' algebriche, ed in
specie delle superficie.

Osserviamo preliminarmente che un sottoinsieme di $n$ punti distinti
$ \{w_1, \dots w_n \} \subset \C$ corrisponde biiettivamente ad un
polinomio monico 
$$ P(z) := ( \prod_{i=1}^n  (z- w_i) ) \in \C [z]_n.$$
Viceversa, un polinomio monico $P(z) \in \C[z]_n$ ha una fattorizzazione come
sopra in $n$ fattori lineari, solo che le $n$ radici sono distinte se e 
solo se il discriminante
$\delta(P)$ (una funzione polinomiale dei coefficienti di $P$) non
si annulla.

\begin{df}
Il {\bf gruppo delle trecce di Artin}, denotato $\BB_n$, 
e' il gruppo fondamentale 
$$ \BB_n : =  \pi_1 ( \C[z]_n - \{ P| \delta (P) = 0\}) $$
dello spazio dei polinomi monici di grado $n$ a radici distinte.
\end{df}

In genere si prende come "punto base" in tale spazio il polinomio 
$ P(z) := ( \prod_{i=1}^n  (z- i) ) \in \C [z]_n$, ovvero, l'insieme
$\{ 1,2,\dots n\}$. Dato un cammino chiuso $\PPP : [0,1] \ra 
( \C[z]_n - \{ P| \delta (P) = 0\})$, si associa ad esso il sottinsieme
di $\R^3 = \C \times \R$ dato da $ \{ (z,t) | \PPP(t) (z) := P_t (z)= 0 \}. $ 

La Figura 3 mostra appunto due realizzazioni della
stessa treccia. Si puo' inoltre osservare, 
come suggerisce la figura, che, poiche' esiste
un sollevamento di $\PPP$ allo spazio $\C^n$ delle $n$-uple di radici di 
un polinomio, si possono cio\'e trovare  funzioni $w_i(t)$ tali che

$ w_i(0) = i$ e $\PPP(t) (z) := P_t (z)= \prod_{i=1}^n ( z - w_i(t))$

ad una treccia corrisponde anche una permutazione $\tau \in \SSS_n$ definita da
$ \tau(i) := w_i(1) \in \{ 1,2, \dots n\}.$

\begin{figure}[htbp]
\begin{center}
\scalebox{1}{\includegraphics{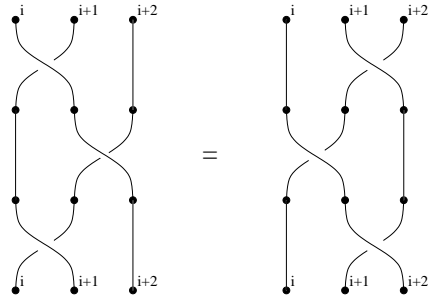}}
\end{center}
\caption{Relazione $aba=bab$ sulle trecce}
\label{figura3}
\end{figure}

Per quanto visivamente suggestiva sia la rappresentazione 
delle trecce data da Artin, una generalizzazione potente e' stata data
da Max Dehn nel suo lavoro \cite{dehn} (rimandiamo il lettore per questo
ad altri punti che ora menzioneremo al libro \cite{birm}).

\begin{df}
Data una varieta' differenziabile $M$, il {\bf gruppo del Dehn} 
(delle classi di applicazioni) di $M$,
denotato con $\MM ap (M)$, e' il gruppo quoziente 
$ \pi_0 (\D iff (M)) = \D iff (M) / \D iff^0 (M)$
delle componenti connesse per archi del gruppo dei diffeomorfismi di $M$.

Il sottogruppo $ \D iff^0 (M)$ si chiama il sottogruppo dei diffeomorfismi
{ \bf isotopi}  all' identita'.

In pratica, quando la varieta' $M$ e' orientata (ad esempio per le varieta'
complesse la struttura complessa induce una orientazione naturale), si considera
tacitamente il sottogruppo $\D iff (M)^+$ dei diffeomorfismi che conservano la orientazione.

In questo caso dunque e' meglio definire $\MM ap (M) =
\pi_0 (\D iff (M)^+) = \D iff (M)^+ / \D iff^0 (M)$.

Nel caso particolare ove $M$ sia una curva complessa compatta di genere $g$,
il Mapping Class group si denota con $\MM ap_g$.
\end{df}

La relazione che sussiste fra gruppo di Artin delle trecce e gruppo del
Dehn e' indicata dal seguente

\begin{teo}
Il gruppo delle trecce $\BB_n$ e' isomorfo al gruppo 
$$  \pi_0 (\MM ap^{\infty} ( \C - \{1, \dots n\}))$$ 

delle componenti connesse
del gruppo $\MM ap^{\infty} ( \C - \{1, \dots n\})$ dei diffeomorfismi che 
sono l' identita' al di fuori del cerchio di centro l' origine
e raggio $2n$.
\end{teo}

Tramite questo isomorfismo si rappresentano i generatori standard di Artin
$\sigma_i$ del gruppo $\BB_n$ ( $ i=1, \dots n-1$) tramite i cosiddetti
half-twists o mezzi-giri:

\begin{df}
Il {\bf mezzo-giro} o {\em half-twist} $\sigma_j$ e' il diffeomorfismo
del piano complesso che e' dato
da una rotazione di 180 gradi sul cerchio di centro $j+1/2$ e raggio $1/2$ ,
e che poi su una circonferenza con lo stesso
centro e raggio $ 1/2 + t/4$ e' dato dall'identita' se $ t \geq 1$, e da
una rotazione
di angolo $180 (1-t)$ gradi se  $ t \leq 1$.
\end{df}
 
Poiche' $\BB_n$ e' visto come sottogruppo di
$\MM ap ( \C - \{1, \dots n\})$ e' ovvio che si ha una azione di
$\BB_n$ sul gruppo libero $\pi_1( \C - \{1, \dots n\})$ che, 
se si prende come punto base il numero complesso $p := -2ni$, ha una base geometrica
$\gamma_1, \dots \gamma_n$ come descritto in precedenza (si noti che
qui i $\gamma_1, \dots \gamma_n$ non soddisfano nessuna relazione,
si ha la relazione  $ \Pi_{h=1}^n \gamma_h = 1 $ solo quando si  considerano
 questi cammini come elementi di 
$\pi_1( \PP^1_{\C} - \{1, \dots n\})$).

La azione cosi' ottenuta si chiama

{\bf Azione di Hurwitz del gruppo delle trecce}, che sostanzialmente 
corrisponde ai cambiamenti di base geometrica (parimenti il gruppo lineare
di uno spazio vettoriale corrisponde ai cambiamenti di base vettoriale)
ed ha la seguente descrizione algebrica:

\begin{itemize}
\item
$\sigma_i (\gamma_i) = \gamma_{i+1}$
\item
$\sigma_i (\gamma_i \gamma_{i+1}) = \gamma_i \gamma_{i+1}$ e quindi
$\sigma_i (\gamma_{i+1}) = \gamma_{i+1} ^{-1}\gamma_i \gamma_{i+1}$
\item 
$\sigma_i (\gamma_j ) = \gamma_j $ per $ j \neq i, i+1$.
\end{itemize}

Si osserva, come e' intuitivo geometricamente, che il prodotto 
$\gamma_1  \gamma_2 \dots    \gamma_n $ e' lasciato invariante
dall' azione del gruppo, e questa osservazione conduce ad una
osservazione assai generale

\begin{df}
Sia $G$ un gruppo, e consideriamo il suo prodotto Cartesiano 
$G^n$. La applicazione che ad una n-upla ordinata $(g_1, g_2, \dots g_n)$
fa corrispondere il
prodotto $g := g_1 g_2 \dots g_n \in G$ ci da' una partizione di $G^n$
i cui sottinsiemi si chiamano {\bf fattorizzazioni} di un elemento 
$ g \in G$. 
Il gruppo delle trecce $\BB_n$ agisce dunque su $G^n$ preservando  la partizione,
e le sue orbite si chiamano {\bf classi di equivalenza di Hurwitz} di fattorizzazioni.

\end{df} 

Per ricordare il concetto di fattorizzazione useremo la notazione
$ g_1 \circ g_2  \circ \dots  \circ g_n $ che distingue quindi la 
fattorizzazione dal suo prodotto $g_1 g_2 \dots g_n$ 
(l' ultima notazione denota invece
l' elemento corrispondente di $G$).

\begin{oss}
Una equivalenza piu' larga sull' insieme delle fattorizzazioni
si ottiene considerando l' equivalenza generata dalla equivalenza di Hurwitz
e dalla {\bf coniugazione simultanea} che, usando la notazione 
 $ a_b : = b^{-1} a b$, corrisponde alla azione di $G$
su $G^n$ che porta  $ g_1 \circ g_2  \circ \dots  \circ g_n $ in
$ (g_1)_b \circ (g_2)_b  \circ \dots  \circ (g_n)_b$.

Si noti che tale azione porta una fattorizzazione di $g$ in una
fattorizzazione del coniugato $g_b$ di $g$.
\end{oss}

La equivalenza di cui sopra gioca un ruolo determinante in varie questioni sulle
curve piane e sulle superficie algebriche, come vedremo nei prossimi
paragrafi. Procediamo intanto a vedere una altra interessante relazione fra il gruppo
delle trecce ed il Mapping class group.

Tale relazione e' strettamente collegata col modello topologico
di una superficie di Riemann di genere $g$ provvisto dalla
curva iperellittica $C_g$ di equazione 
$$ w^2 = \prod_{i=1}^{2g + 2}  ( z - i)  $$ 
(vedi la Figura 4 che descrive una curva iperellittica di genere $g=2$).

\begin{figure}[htbp]
\begin{center}
\scalebox{1}{\includegraphics{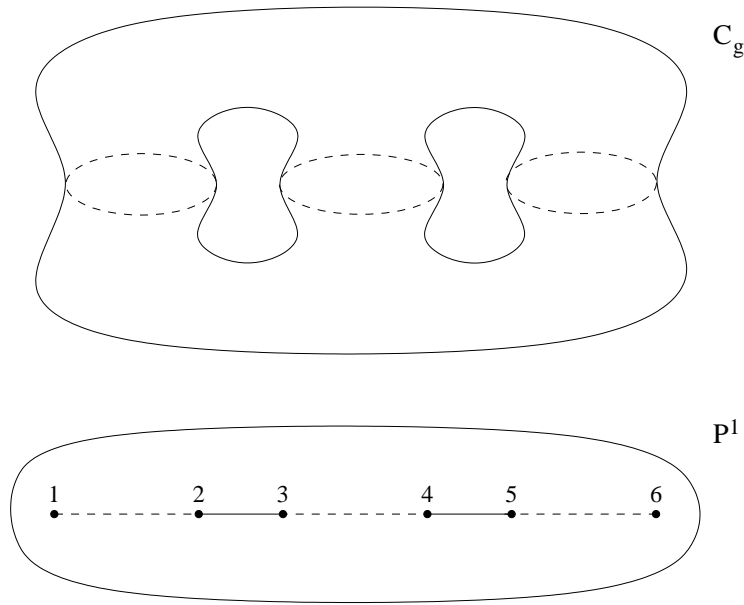}}
\end{center}
\caption{Curve iperellittica di genere $2$}
\label{figura4}
\end{figure}

Si osservi che $C_g$ e' la compattificazione naturale dei limiti all'infinito
({\bf ends} in inglese) del rivestimento  doppio $Y$ di $( \PP^1 - \{1, \dots 2g+2\})$  
dato dalla controimmagine di $( \PP^1 - \{1, \dots 2g+2\})$ in $C_g$
(in tale compattificazione  si aggiunge a $Y$  il $ lim_{ K \subset \subset Y}
\pi_0 (Y - K)$).

E' chiaro dunque che ogni omeomorfismo di $( \PP^1 - \{1, \dots 2g+2\})$
che lascia stabile il sottogruppo 
associato al rivestimento $Y$ si lascia sollevare ad un omeomorfismo di $Y$,
e quindi anche della sua compattificazione naturale $C_g$.

Un tale sollevamento non e' unico, perche' possiamo sempre comporre con un
automorfismo del rivestimento.

Abbiamo quindi una successione esatta (piu' precisamente una estensione centrale)

$$ 1 \ra \Z /2 =  <H> \ra \MM ap_g^h  \ra \MM ap_{0,2g+2} \ra 1 $$

ove
\begin{itemize}
\item
$H$ e' la involuzione iperellittica $ w \ra - w$,
 che e' l' unico automorfismo non banale del rivestimento
\item
 $\MM ap_{0,2g+2}$ e' il  gruppo del Dehn di $( \PP^1 - \{1, \dots 2g+2\})$
\item
$\MM ap_g^h$ si chiama il sottogruppo iperellittico del
mapping class group $\MM ap_g$, e consiste appunto di tutti i possibili sollevamenti.

Se $ g \geq 3$, e' un sottogruppo proprio di $\MM ap_g$.

\end{itemize}

E' interessante osservare che, mentre il gruppo delle trecce di Artin
  $\BB_{2g+2}$ ha la seguente presentazione:
$$ < \sigma_1,  \dots \sigma_{2g+1} |  
 \sigma_i  \sigma_j =  \sigma_j
 \sigma_i  \ per  |i-j| \geq 2, \ \sigma_i  \sigma_{i+1}  \sigma_i =
\sigma_{i+1}  \sigma_i \sigma_{i+1}>, $$ 

il gruppo del Dehn di $( \PP^1 - \{1, \dots 2g+2\})$ ha la presentazione $\MM ap_{0,2g+2}$ :
$$ < \sigma_1,  \dots \sigma_{2g+1} |   \sigma_1  \dots \sigma_{2g+1}
 \sigma_{2g+1} \dots  \sigma_1 = 1 , (\sigma_1  \dots \sigma_{2g+1})^{2g+2} =1,$$ $$
 \sigma_i  \sigma_j =  \sigma_j
 \sigma_i  \ per  |i-j| \geq 2, \ \sigma_i  \sigma_{i+1}  \sigma_i =
\sigma_{i+1}  \sigma_i \sigma_{i+1}>, $$ 
infine il mapping class group iperellittico  $ \MM ap_g^h $ ha la presentazione:

$$ < \xi_1,  \dots \xi_{2g+1}, H |   \xi_1  \dots \xi_{2g+1}
  \xi_{2g+1} \dots  \xi_1 = H , H^2 = 1,(\xi_1  \dots \xi_{2g+1})^{2g+2} =1, $$ $$
  H \xi_i = \xi_i H \ \forall i,  \xi_i  \xi_j =  \xi_j
 \xi_i  \ per  |i-j| \geq 2, \ \xi_i  \xi_{i+1}  \xi_i =
\xi_{i+1}  \xi_i \xi_{i+1} >. $$

Ma conviene forse lasciare da parte queste importanti formule e ritornare alla
geometria.
Osserviamo che $\sigma_j$ da' un omeomorfismo del cerchio $U$ di centro $ j + 1/2$ 
e raggio 3/4, inoltre scambia tra di loro i due punti $j, j+1$. 

Quindi ci sono due sollevamenti di $\sigma_j$ ad omeomorfismi della
controimmagine $V$ di $U$ in $C_g$: si definisce allora $\xi_j$ come quello dei due
sollevamenti che agisce come l'identita' sul bordo di $V$ (che e' una
unione di due circonference,
vedi figura 5). 

$\xi_j$  si chiama {\bf GIRO del Dehn, o Dehn twist} e 
corrisponde geometricamente al diffeomorfismo di un tronco di cilindro che
e' l' identita' sul bordo, una rotazione di 180 gradi all' equatore, e su ogni
parallelo  ad altezza $t$ una rotazione di $ t \  360 $ gradi ( qui $t \in [0,1]$).

\begin{figure}[htbp]
\begin{center}
\scalebox{1}{\includegraphics{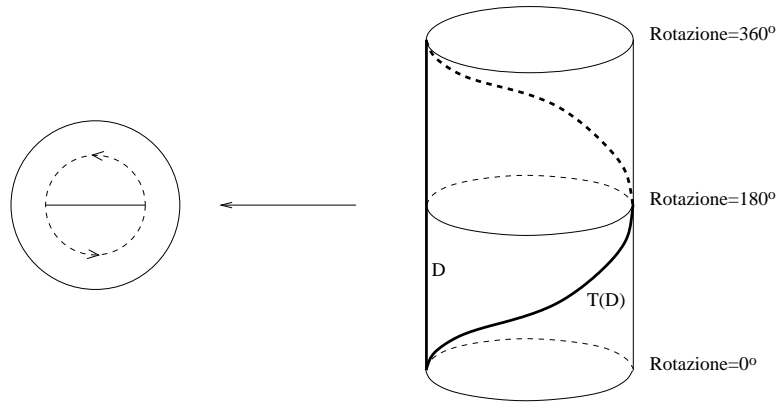}}
\end{center}
\caption{A sinistra il mezzo giro, a destra il suo sollevamento: 
il Dehn-Twist-$T$ e la sua azione sul segmento $D$}
\label{figura5}
\end{figure}

In realta' si puo' definire analogamente un Dehn twist per ogni curva orientata
in $C_g$ (cioe' una sottovarieta' diffeomorfa ad $S^1$), ed il risultato
fondamentale del Dehn (\cite{dehn}) e' stato il

\begin{teo}
Il Mapping Class group $\MM ap_g$ e' generato dai Dehn twists.
\end{teo}

Delle presentazioni esplicite di $\MM ap_g$ sono in seguito state date
da Hatcher e Thurston (\cite{h-t}) e con un miglioramento del metodo
la presentazione piu' semplice e' dovuta a Wajnryb (\cite{waj}, vedi anche
\cite{waj2}).

Vedremo nel prossimo paragrafo come i twists del Dehn sono collegati
con la teoria delle fibrazioni di Lefschetz.

\section {Fibrazioni di Lefschetz.}

Il metodo introdotto da Lefschetz per l'indagine delle proprieta' topologiche
delle varieta' algebriche e' l' analogo topologico dei metodi di
proiezione e sezione della scuola italiana classica di geometria algebrica.

Nel caso delle varieta'  compatte di dimensione reale $4$ i metodi
della teoria di Morse e della semplificazione del cobordismo
avevano trovato ostacoli estremamente ardui da superare, e solo nel 1982
M. Freedman (\cite{free}) , usando nuove idee per dimostrare la trivialita' 
topologica di certi manici considerati da A. Casson, ha ottenuto 
una classificazione delle 4-varieta' topologiche compatte semplicemente connesse
che ha come corollario

\begin{teo}
Sia $S$ una superficie complessa compatta e semplicemente connessa,
e sia $r$ l' indice di divisibilita' della classe canonica 
$c_1( K_X) \in H^2 (X, \Z)$.

$S$ si dice PARI se $ r \equiv 0(mod 2)$, DISPARI altrimenti.
Allora
\begin{itemize}
\item
 Se $S$ e' PARI (questo vale se e solo se la forma di intersezione su  $H^2 (X, \Z)$
e' pari) allora  $S$ e' una somma connessa di copie di
 $ \PP ^1_{\C} \times \PP ^1_{\C}$ e di una superficie $K3$ 
se la segnatura  della forma di intersezione e' negativa, e di copie di
$ \PP ^1_{\C} \times \PP ^1_{\C}$ e di una superficie $K3$ con
orientazione rovesciata nel caso che la detta segnatura sia positiva.
\item
$S$ e' dispari: allora $S$ e' una somma connessa di copie di
$\PP ^2_{\C} $ and ${\PP ^2_{\C}}^{opp} .$
\end{itemize}
\end{teo}

\begin{oss}
${\PP ^2_{\C}}^{opp} $ e' la  varieta' 
${\PP ^2_{\C}}$con orientazione rovesciata.

Invece, una superficie K3 e' una superficie  $S$ orientatamente 
diffeomorfa ad una superficie
$X$ non singolare di grado 4 in $\PP^3_{\C}$,
ad esempio 
$$ X = \{ (x_0,x_1,x_2,x_3) 
\in \PP^3_{\C} | x_0^4 + x_1^4 + x_2^4 + x_3^4 = 0\}.$$
(per un teorema di Kodaira, cf. \cite{kod},  $S$ e' anche equivalente per deformazione a tale $X$).
\end{oss} 

La teoria di Donaldson ha reso chiaro in seguito  (\cite{don1}, \cite{don2},
\cite{don3},\cite{don4}) quanto drasticamente omeomorfismo e 
diffeomorfismo differiscano in dimensione 4 (a differenza delle altre dimensioni),
e specialmente per le superficie algebriche.

La citata congettura di
Friedman and Morgan  (Def = Diff) fu motivata dai risultati positivi ottenuti 
tramite la teoria di gauge e gli invarianti di Donaldson.

In seguito, la teoria di Seiberg e Witten ha mostrato con metodi piu'
semplici che ogni diffeomorfismo tra superficie minimali (ad esempio, nel caso PARI
ogni superficie e' necessariamente minimale)
porta 
 $c_1( K_X)$ su $c_1( K_{X'})$ o su  $- c_1( K_{X'})$
(cf. \cite{witten} o \cite{mor}).

Ne consegue quindi che l' invariante $r$, definito sopra tramite la divisibilita'
della classe canonica, e' un invariante differenziale: ma solo la sua parita' e'
un invariante topologico; quindi e' facile trovare esempi di superficie algebriche
semplicemente connesse che siano fra loro omeomorfe ma non diffeomorfe
(vedi gli esempi di \cite{cat3}). 

Una questione dunque assai difficile e' quella di decidere, per esempio anche nel caso
semplicemente connesso,  se due superficie algebriche, non deformazione
l'una dell' altra, siano effettivamente tra loro diffeomorfe.

Per questo scopo la teoria delle fibrazioni di Lefschetz offre un criterio di
determinazione positiva,  
dovuto a Kas (\cite{kas}), e che andiamo ad illustrare. Ci soffermeremo per brevita'
ad illustrare solo il caso di fibrazioni di Lefschetz in dimensione reale 4.

\begin{df}
Sia $M$ una 4-varieta' differenziabile compatta.
Possiamo piu' generalmente
supporre che su $M$ sia data una struttura simplettica, cioe' una 
2-forma differenziale $\omega$ tale che $ d (\omega) = 0$
e tale che $\omega$ dia una forma alternante non degenere in ogni puto di $M$.

Una fibrazione di Lefschetz e' una applicazione differenziabile 
 $ f : M \ra \PP^1_{\C}$
che ha differenziale di rango massimo eccetto per
un insieme finito di punti critici $p_1, \dots p_m $
che hanno valori critici distinti
$b_1, \dots b_m \in 
\PP^1_{\C}$, ed ha la proprieta' che intorno al punto $p_i$
esistono coordiante complesse $(x,y) \in \C^2$  tale che localmente 
$f = x^2 - y^2 + costante $ 
(nel caso simplettico, nelle date coordinate  la forma
$\omega$ deve corrispondere alla struttura simplettica standard di
 $\C^2$).
\end{df}

Una definizione simile si puo' dare se $M$ ha bordo, 
rimpiazzando $\PP^1_{\C}$ con un disco $ D \subset \C$.

Un importante teorema di Donaldson (\cite{don6}) assicura che 
per le varieta' simplettiche vale, come per il caso di
varieta' algebriche proiettive, la esistenza di un fascio di Lefschetz,
cioe' una fibrazione di Lefschetz $ f : M' \ra \PP^1_{\C}$
su uno scoppiamento simplettico $M'$ di $M$
(cf.\cite{MD-S}).

Vediamo ora come  i Dehn twists intervengono nel contesto
delle fibrazioni di Lefschetz.

Consideriamo infatti le coordinate locali speciali intorno ad un punto 
critico, e supponiamo per semplicita' (a meno di traslazione) 
che il valore critico sia uguale a $0$.

Allora, al variare del numero complesso $ t$ in un disco
di raggio $\epsilon$, la fibra $f^{-1}(t)$ e' localmente descritta
dalla equazione  $ y^2 = x^2 - t $, e si ha dunque un rivestimento doppio
diramato sulle due radici $\pm \sqrt t$ del polinomio $x^2 - t $.

Consideriamo la immagine inversa della circonferenza di centro l' origine e raggio
$\epsilon$:
$ f^{-1}( \{ t| t =\epsilon e^{2 \pi i \theta}\})$ e' localmente omeomorfa
al quoziente del prodotto Cartesiano $ f^{-1}(\epsilon) \times [0,1]$ 
per la relazione di equivalenza che identifica $ f^{-1}( \epsilon)
\times \{ 0 \}$ con $ f^{-1}( \epsilon)
\times \{ 1 \}$  tramite il diffeomorfismo $ T$ dato dal Dehn 
twist sul ciclo evanescente immagine inversa del segmento con estremi
$\pm \sqrt \epsilon$ (il ciclo si chiama evanescente 
perche' chiaramente,
quando $t$ tende a zero, il segmento di estremi $\pm \sqrt t$ tende a
svanire in  un punto).

In altre parole, la {\bf monodromia locale}, cioe' sopra la circonferenza
 $\{ t| t =\epsilon e^{2 \pi i \theta}\})$ e' dato dal sollevamento 
del mezzo giro sul segmento
di estremita' $\pm \sqrt \epsilon$ (vedi figura 5).

Quel che succede globalmente e' che una fibrazione di Lefschetz di curve di genere 
$g$ e con valori critici
$b_1, \dots b_m \in 
\PP^1_{\C}$, una volta che si e' scelta una quasi-base geometrica 
$\gamma_1 , \gamma_2 ,  \dots \gamma_m $ di
$ \pi_1 ( \PP^1_{\C} - \{ b_1, \dots b_m \}) $,
 determina una fattorizzazione della identita'  
$$\tau_1 \circ \tau_2 \circ \dots \tau_m = Id  $$
come prodotto di Dehn twists nel  Mapping class group
$  \MM ap_g .$ 

Possiamo ora enunciare l' importante risultato ottenuto da Kas nel 1980 (\cite{kas}) 
(influenzato dalla rivisitazione da parte di Andreotti-Frenkel e 
Moisezon dei risultati di Lefschetz)

\begin{teo}{\bf (Kas)}
Due fibrazioni di Lefschetz  $(M, f)$,
$(M', f')$
si dicono equivalenti se esistono due diffeomorfismi
$u : M \ra M', v : \PP^1 \ra \PP^1$ tali che $ f\ \circ u = v \circ f$.

Orbene, la classe di equivalenza di una fibrazione di 
Lefschetz $(M, f)$ e' completamente determinata dalla classe di equivalenza
della sua fattorizzazione della identita' nel
Mapping class group (per la relazione di equivalenza generata da 
equivalenza di Hurwitz e da coniugazione simultanea).

\end{teo}

Un risultato simile vale per una fibrazione di Lefschetz
sopra un disco $ D \subset \C$:
otteniamo una fattorizzazione
$$\tau_1 \circ \tau_2 \circ \dots \tau_m = \phi $$
della monodromia $\phi$ della fibrazione sulla circonferenza bordo di $D$.

Inoltre, la fibrazione ammette una struttura simplettica se e solo se
ogni $\tau_i$ e' un Dehn twist positivamente orientato.

Siamo ora in grado di ricucire il filo con gli esempi di chirurgia
mostrati nel primo paragrafo.

Infatti, date due fibrazioni di Lefschetz su $\PP^1_{\C}$, 
consideriamo come sottovarieta'
$N_i$ una fibra liscia della fibrazione, dimodoche' il suo fibrato
normale e' banale : allora la chirurgia associata alle trivializzazioni
da' un risultato che si chiama {\bf Somma sulla fibra} ({\bf fiber sum} in inglese) 
delle due fibrazioni di Lefschetz., e dipende dal
diffeomorfismo scelto tra $N_1$ ed $N_2$ (cf. \cite{g-s},
 Def. 7.1.11, e Teorema 10.2.1.).

Questa operazione si traduce, in vista del risultato di Kas
sopra citato, nella composizione "coniugata" di due fattorizzazioni

\begin{df}
Siano $\tau_1 \circ \tau_2 \circ \dots \circ \tau_m = \phi $ e
$\tau'_1 \circ \tau'_2 \circ \dots \circ \tau'_r = \phi' $ 
due fattorizzazioni: allora la loro composizione coniugata per
$\psi$ e' la fattorizzazione
$$\tau_1 \circ \tau_2 \circ \dots \tau_m \circ
 (\tau'_1 )_{\psi}\circ (\tau'_2)_{\psi} \circ \dots \circ 
(\tau'_r)_{\psi}= \phi (\phi')_{\psi}. $$
Se ${\psi}$ commuta con $\phi'$, otteniamo una fattorizzazione di $\phi \phi'$.

Un caso particolare e' quello in cui $\phi, \phi'$ sono banali,
caso che corrisponde alle fibrazioni di Lefschetz su $\PP^1$.

\end{df} 

\section { ABC e Def $\neq$ Diff nel caso semplicemente connesso.}

Abbiamo osservato che, se due varieta' complesse $X, X'$  sono 
equivalenti per deformazione, allora esiste un diffeomorfismo 
che porta la classe canonica  $c_1( K_X) \in H^2 (X, \Z)$ di $X$
su quella di $X'$. D' altra parte per il risultato di Seiberg-Witten
  un diffeomorfismo
porta la classe canonica di una superficie minimale $S$ su
$\pm c_1( K_{S'})$. Quindi, se diamo almeno tre superficie  fra
loro (orientatamente) diffeomorfe,
ne troviamo sicuramente due tra cui esiste un diffeomorfismo con la proprieta'
di portare la classe canonica sulla classe canonica.

Il seguente teorema, ottenuto in \cite{c-w}, da' il contresempio cercato
alla congettura di Friedman Morgan anche nel caso semplicemente connesso.

\begin{teo}
Per ogni intero positivo $h$, esistono
$h $ superficie  $S_1, \dots S_h$ semplicemente connesse,
diffeomorfe l' una con l'altra, ma tali che due
di queste superficie 
non sono mai equivalenti per deformazione.
\end{teo}

Gli esempi considerati sono le superficie $(a,b,c)$ ottenute come
compattificazioni minime delle superficie affini di equazioni

\begin{eqnarray}  z^2 &=&
 f(x,y)\\
 w^2 &=&
 g(x,y) \nonumber ,\end{eqnarray}

dove f e g sono polinomi di bigradi rispettivi $(2a,2b), (2c, 2b)$.

Per ottenere le compattificazioni si considerano le estensioni di questi
polinomi a sezioni di fasci invertibili su $\PP^1\times \PP^1$:
 $ f \in H^0({\PP^1\times \PP^1}, {\hol}_{\PP^1\times \PP^1}(2a,2b)) $ e  
$ g \in H^0({\PP^1\times \PP^1}, {\hol}_{\PP^1\times \PP^1}(2c,2b)), $
e se ne prendono le radici quadrate dentro lo spazio totale 
della somma diretta dei due fibrati lineari associati a 
${\hol}_{\PP^1\times \PP^1}(a,b), {\hol}_{\PP^1\times \PP^1}(c,b)$.

La superficie non singolare e compatta $S$  ottenuta e' un rivestimento "bidoppio"
(cioe' Galois con gruppo $(\Z/2)^2$ )
 di ${\Bbb P}^1
\times  {\Bbb P}^1  $ semplice di 
 tipo  ((2a, 2b),(2c,2b), e viene denotata per semplicita' come "superficie $(a,b,c)$".

La compattificazione $S$ e' non singolare purche' le due curve
$f=0$ e $g=0$ siano lisce e si intersechino trasversalmente.

Il teorema fondamentale sopra citato segue  piu' precisamente 
dai due teoremi seguenti:

\begin{teo}
Siano  $S$ una superficie $(a,b,c)$, $S'$ una superficie 
$(a+k,b,c-k)$.

Se assumiamo 

\begin{itemize}
\item
 (I) $a,b,c, k$
sono interi pari positivi con $ a, b, c-k \geq 4$
\item
(II)  $ a \geq 2c + 1$,
\item
(III) $ b \geq c+2$ .
\end{itemize}
Allora  $S$ e  $S'$ non sono equivalenti per deformazione
\end{teo}

\begin{teo}
Siano $S$, $S'$ una superficie $(a,b,c)$, rispettivamente  una superficie 
$(a+1,b,c-1)$. , e si assuma 
 che $a,b,c-1$
siano interi $ a, b, c-1 \geq 2$.

Allora   $S$ e  $S'$ sono diffeomorfe.
\end{teo}

Si noti che le superficie in questione sono semplicemente connesse,
come mostrato in (cf.  \cite [Proposition 2.7]{cat1}).

Il primo di questi teoremi usa tecniche e risultati diversi, sviluppati
in una serie di lavori (\cite{cat1}, \cite{cat2},\cite{cat3},
\cite{man1},\cite{man3}), e che appartengono alla geometria algebrica ed 
alla teoria delle singolarita' : teoria delle piccole deformazione alla Kuranishi,
degenerazioni  normali di superficie lisce, singolarita' quozienti
di punti doppi razionali.

Un ingrediente estremamente semplice ma di fondamentale importanza e' il
concetto di
{\bf  deformazioni naturali
di un rivestimento bidoppio } (\cite{cat1}, def. 2.8 , page 494 )
che sono parametrizzate da una quaterna di polinomi
$f, g, \phi, \psi$ a cui corrispondono le equazioni
\begin{eqnarray}  z^2 &=&
 f(x,y) + w \phi(x,y)\\
 w^2 &=&
 g(x,y) + z \psi(x,y) \nonumber ,\end{eqnarray}

($f,g$ sono come prima , mentre $\phi$ ha bigrado (2a-c,2b-d),
$\psi$ ha
bigrado (2c-a,2d-b)).

Le deformazioni naturali danno tutte le piccole deformazioni, sotto 
ipotesi opportune (base rigida, e curve di diramazione abbastanza positive).

Inoltre, poiche' $ a \geq 2c+1$, segue che  $ \psi \equiv 0$,
 ed allora ogni piccola deformazione conserva la struttura di 
{\em rivestimento doppio iterato } (\cite {man3}): cioe' ciascuna
superficie possiede l' automorfismo $ z \ra -z$, ed il quoziente per questo
automorfismo possiede l' automorfismo $ w \ra -w$. 

Il punto finale e' di mostrare che questa struttura passa convenientemente
al limite, in modo che si ottiene un inseme non solo aperto, ma
anche chiuso dello spazio dei moduli.

Sugli altri metodi assai "complessi" utilizzati non mi voglio qui soffermare, 
sia per mancanza di tempo, sia
 perche' scopo di questo articolo
era di illustrare i metodi di topologia differenziale usati per il secondo
di questi due teoremi.

Le idea chiave usate per quest' ultimo sono che: 

1) entrambe le superficie $S$ ed $S'$ posseggono
una applicazione olomorfa su 
 $\PP^1_{\C}$ data dalla composizione del rivestimento bidoppio 
colla proiezione sulla prima coordinata $x$, ed una piccola perturbazione 
di questa applicazione le realizza come fibrazioni di
 Lefschetz simplettiche (cf. \cite{don7}, \cite{g-s})

2) le rispettive fibrazioni, proprio per la scelta accorta dei gradi delle
curve (la prima curva $\{ f=0 \}$ perde un bigrado (2,0), e la seconda  
$\{ g=0 \}$
lo acquista) sono somma sulla fibra della stessa coppia di fibrazioni
di Lefschetz sopra il disco

3) per dimostrare  che le due  somme sulla fibra sono equivalenti e' sufficiente 
mostrare che, una volta che la prima delle due sia rappresentata come
composizione di due fattorizzazioni, e la seconda come composizione 
coniugata per $\psi$ delle stesse fattorizzazioni, 

(**) il diffeomorfismo 
$\psi$ sta nel sottogruppo del Mapping Class Group della fibra
generato dai Dehn twists che compaiono 
nella prima fattorizzazione (lemma di Auroux, \cite{aur}, cf. anche
\cite{kas}).

4) La figura 6 mostra la curva $C$ fibra della fibrazione nel caso $2b=6$:  essa
e' un rivestimento bidoppio di $\PP^1_{\C}$, che possiamo supporre 
 dato dalle equazioni

\begin{itemize}
\item
$ z ^2 = F(y) , w  ^2 = F (-y)$\\
dove le radici di $F$ sono gli interi positivi $1, \dots 2b$,
\item
si vede inoltre che la monodromia della fibrazione sul bordo del disco e' banale 

\item
 e che la applicazione $\psi$ e' data dal diffeomorfismo
di ordine due

%\item
$ (*) \  y \ra -y , z \ra w , w \ra z$

\end{itemize}
che nella figura 6 e' rappresentato da una rotazione di 180 gradi
intorno ad un asse inclinato verso Nord Ovest.

La figura mostra una evidente simmetria diedrale, ove l' automorfismo
di ordine $4$ e' dato da $ (*) \  y \ra -y , z \ra - w , w \ra z$.

Inoltre, fra  i Dehn twists che compaiono nella fattorizzazione ci sono
quelli corrispondenti alle curve immagini inverse dei segmenti
fra due interi consecutivi (vedi figura 6). 

\begin{figure}[htbp]
\begin{center}
\scalebox{1}{\includegraphics{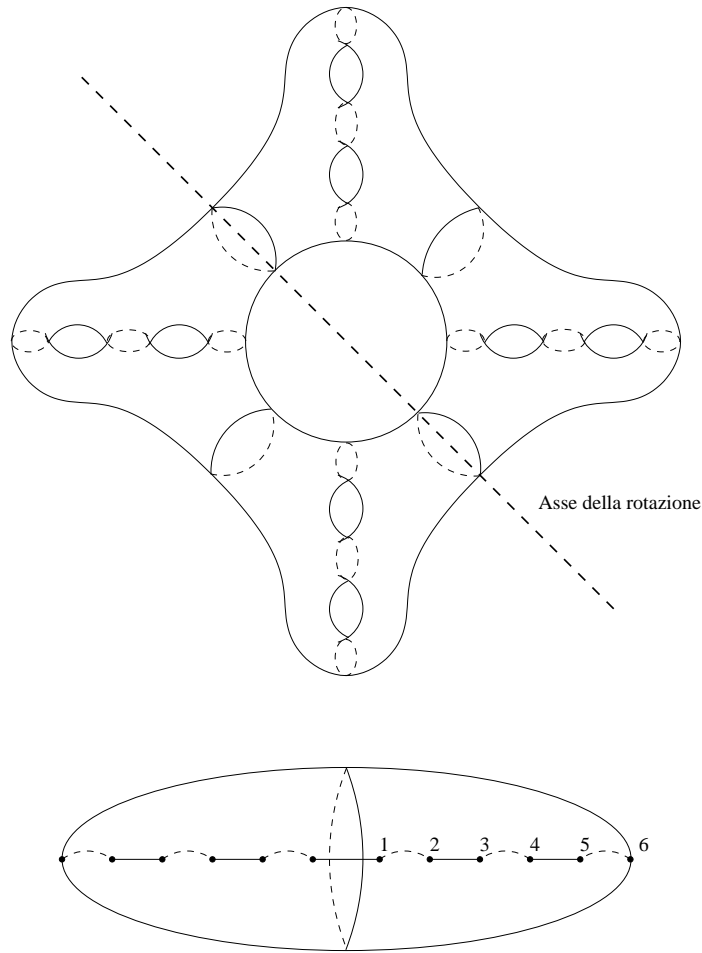}}
\end{center}
\caption{La Curva Complessa $C$ a Simmetria Diedrale}
\label{figura6}
\end{figure}

Queste curve si possono organizzare sulla curva $C$ in sei catene
(non tra loro disgiunte), ed infine si e' ridotti a mostrare
che la classe di isotopia di $\psi$ e' la stessa del prodotto 
di sei elementi di Coxeter associati a tali catene.

Ricordiamo che gli elementi di Coxeter di una catena sono 
prodotti del tipo 

$\Delta= (T_{\alpha_1})(T_{\alpha_2}T_{\alpha_1})\dots
(T_{\alpha_{n-1}}T_{\alpha_{n-2}}\dots T_{\alpha_1})(T_{\alpha_n}
T_{\alpha_{n-1}}\dots T_{\alpha_1})$ dei Dehn twists 
associati alle curve della catena.

Per dimostrare infine che tale prodotto $\psi'$ di elementi di Coxeter e $\psi$
sono tra loro isotopi, si osserva che rimuovendo le curve sopracitate dalla
curva complessa $C$ si ottengono 4 componenti connesse diffeomorfe a dei cerchi.
Basta allora verificare (per un risultato di Epstein) che $\psi'$  e $\psi$
portano ciascuna di tali curve su una coppia di curve isotope:
l' ultimo passo richiede una serie di lunghe  anche se facili 
 verifiche, per cui sono importanti delle figure esplicite.

\section {Varieta' simplettiche e monodromia delle trecce.}

Nel lavoro (\cite{cat6}) e' stato osservato che una superficie di tipo
generale possiede una struttura simplettica canonica, 
unica a meno di simplettomorfismo, tale che la classe della forma simplettica
 $\omega$ in coomologia di De Rham sia l'immagine della  classe canonica 
  $c_1( K_X) $ in $ H^2 (X, \R)$.

Nel migliore dei mondi possibili, assumendo cioe' che il divisore
canonico sia ampio (ed in ogni dimensione), si prende semplicemente il pull-back
della forma di Fubini Study tramite una immersione $m$-pluricanonica in uno
spazio proiettivo, e poi la si
divide per $m$. Il lemma di Moser assicura che la struttura cosi' ottenuta e' unica.

Una domanda naturale e' se le superficie di tipo (a,b,c) fra loro diffeomorfe
(cioe' quelle con gli stessi valori di (a+c) e b) sono anche tra loro
simplettomorfe se munite della struttura simplettica canonica.

Un modo di dimostrarlo analizzando la
monodromia delle trecce della curva di diramazione del rivestimento
quadruplo "perturbato" (corrispondente alla fibrazione di Lefschetz),
e mostrando che  la involuzione $\iota$ su $\PP^1_{/C}$
tale che $\iota (y) = -y$ e' un prodotto delle trecce che capitano nella
fattorizzazione
si e' rivelata assai piu' difficile della corrispondente analisi fatta nel
Mapping class group (il problema e' che si ha un  omomorfismo non iniettivo da 
un sottogruppo " $\SSS_4$-colorato" , cf. \cite{c-w},
del gruppo  delle trecce dentro
il Mapping class group).

Ci si puo' chiedere piu' generalmente se, 
considerando fibrazioni olomorfe $ f : S \ra B$ su una curva $B$, 
con fibre curve di genere $g \geq 2$, che ammettano una perturbazione
in una fibrazione di Lefschetz simplettica con tutte le fibre singolari 
irriducibili (cioe', tali che il nodo 
non disconnette la fibra),
allora la equivalenza come fibrazioni di Lefschetz differenziabili e'
la stessa della equivalenza di fibrazioni di Lefschetz simplettiche.

L' i dea e' di pensare $f$ come data da un morfismo classificante 
$\phi$ da $B$ nello spazio dei moduli compattificato
$\overline{\MM_g}$ delle curve di genere $g$.

Il tipo differenziabile della fibrazione perturbata $f'$ e' 
dato dalla fattorizzazione 
 della monodromia nel mapping class group $\MM ap_g$ . D' altra parte
$\MM ap_g$ agisce in modo propriamente discontinuo
sul dominio di Teichm\"uller $\T_g$ con quoziente
$\MM_g$. La classe della fattorizzazione dovrebbe dare la classe di
isotopia della applicazione perturbata $\phi'$.

Da qui dovrebbe seguire la esistenza ed unicita' di una struttura
simplettica naturale, associata alla classe canonica relativa
$K_{S|B}$, in qualche modo collegata alla struttura 
simplettica naturale introdotta da Gompf 
(\cite{gompf1}) per fibrazioni di Lefschetz.

Si noti che , alla luce del citato risultato di Donaldson che associa
ad una 4- varieta' simplettica un fascio di Lefschetz naturale,
viene fuori che la chiave della classificazione delle 
4- varieta' simplettiche sta nella comprensione delle
classi di equivalenza di fattorizzazioni nel Mapping class group.

Il lavoro \cite{c-w} da' dei risultati positivi, nel senso che
offre la possibilita', in diversi casi importanti, di dimostrare
la equivalenza di due fattorizzazioni.

Ma sarebbe anche interessante dare dei criteri effettivamente
funzionanti e non banali per determinare delle ostruzioni alla
equivalenza di due fattorizzazioni.

Nel prossimo paragrafo vedremo una proposta  fatta recentemente
nel caso analogo
dei rivestimenti generici del piano $\PP^2_{\C}$.

\section{Monodromia delle trecce ed il problema di Chisini.}

Sia $B \subset \PP^2_{\C}$ una curva  algebrica piana di grado $d$,
e prendiamo un punto $P$ generico che non giace su $B$.

Il fascio delle rette $L_t$ passanti per $P$ determina una famiglia ad un parametro
di $d$-uple di punti di $\C \cong L_t -\{P\}$, cioe' $ L_t \cap B$.
Si ottiene quindi una fattorizzazione nel gruppo delle trecce $\BB_d$,
con prodotto questa volta $(\Delta^2)^d $, ove $(\Delta^2) =
(\sigma_{d-1}\sigma_{d-2} \dots \sigma_{1})^d$ e' il generatore
del centro del gruppo delle trecce.

La classe di equivalenza della fattorizzazione (per equivalenza di Hurwitz 
e coniugazione simultanea) non dipende dal punto $P$, se questo e' generico,
e non cambia se $B$ varia in una famiglia equisingolare di curve.

Ad esempio, se $B$ e' liscia, si ha una classe di fattorizzazione data
da mezzi giri. 

Il caso che piu' interessava a Chisini (cf. e.g. \cite{chis1}, \cite{chis2})
era quello delle curve 
{\bf cuspidali}, cioe' che hanno
come singolarita' solo nodi o cuspidi, anche perche' queste sono le
singolarita' del luogo di diramazione di un rivestimento generico 
$ f : S \ra \PP^2_{\C}$, tale cioe' che gli stabilizzatori locali siano
solo $\Z/2 = \SSS_2$ oppure $\SSS_3$, e tali che i punti critici che hanno la stessa
immagine sono al piu' una coppia di ramificazione semplice dove i due rami 
hanno immagine trasversale (e quindi la monodromia locale e' $\SSS_2 \times \SSS_2$).

In questo caso si ha una fattorizzazione 
{\bf di tipo cuspidale}, cioe' in cui tutti i fattori sono potenze
di un mezzo giro con esponente $1,2,3$ a seconda che si abbia una retta
$L_t$ tangente in un punto liscio, una retta per un nodo, oppure una
retta per una cuspide.

I problemi fondamentali relativi alle trecce algebriche 
posti da Chisini sono stati i seguenti

{\bf Congettura di Chisini, cf. \cite{chis1} : }{\em Dati due  rivestimenti generici 
$ f : S \ra \PP^2_{\C}$, $ f' : S' \ra \PP^2_{\C}$ di grado 
maggiore od uguale a $5$, e con la stessa curva di diramazione $B$,
si puo' concludere che i due rivestimenti sono equivalenti?} 

Chisini stesso osserva in \cite{chis1} 
come la condizione che il grado sia almeno 5 sia necessario,
e questo e' collegato col fatto che le tre trasposizioni
$(1,2)(2,3)(3,1) \in \SSS_3$ e le tre trasposizioni
$(1,2)(1,3)(1,4) \in \SSS_4$ soddisfano alle stesse relazioni di commutazione
(e Chisini produce effettivamente due rivestimenti di gradi rispettivi tre
e quattro, generici e con la stessa curva di diramazione).

La congettura e' ancora aperta, ma e' stata sostanzialmente dimostrata
da Viktor Kulikov (\cite{kul}
(con un astuto miglioramento tecnico dovuto a Nemiroski, \cite{nem}):
infatti e' vera se il grado di ciascun rivestimento e' almeno $12$.

 Una risposta negativa, ad opera di Boris Moisezon (\cite{moi94}) ha
avuto invece il seguente

{\bf Problema di Chisini, cf. \cite{chis2} : }{\em Date due  fattorizzazioni 
di tipo cuspidale, rigenerabili alla fattorizzazione
di una curva algebrica piana non singolare, esiste una curva cuspidale
che realizza tale fattorizzazione?} 

Rigenerabile alla fattorizzazione
di una curva algebrica piana non singolare
significa che esiste una fattorizzazione nella classe tale che,
se si rimpiazza ogni fattore del tipo $\sigma^i$ ($i=2,3$), ove
$\sigma$ e' un mezzo giro, con gli $i$ fattori 
corispondenti (ad esempio si rimpiazza $\sigma^3$ con
$\sigma \circ \sigma \circ \sigma$), si ottiene 
la classe di equivalenza 
della fattorizzazione
di una curva algebrica piana non singolare.

Si noti che questa condizione e' necessaria, perche' 
ogni curva piana ha una piccola deformazione che 
e' una curva non singolare.

D' altra parte le famiglie di curve cuspidali di grado $d$ fissato
sono un insieme algebrico, che ha un numero finito di
componenti connesse: ne segue che le classi di
fattorizzazioni in $\BB_d$ provenienti da curve cuspidali
sono solo in un numero finito.

Moisezon mostra invece l'esistenza di infinite fatorizzazioni 
di tipo cuspidale non equivalenti,
poiche' osserva che il gruppo fondamentale del complementare di $B$,
$\pi_1(\PP^2_{\C} - B)$ e' un invariante della fattorizazione.

Come congetturato da Moisezon, il punto e' che una fattorizzazione di tipo
cuspidale, 
ed una monodromia generica a valori in $\SSS_n$ 
 da' luogo ad un rivestimento $ M \ra \PP^2_{\C} $ dove
$M$ e' una 4-varieta' simplettica.

Estendendo le tecniche sviluppate da Donaldson per dimostrare la esistenza di
fibrazioni di Lefschetz, Auroux e Katzarkov (\cite{a-k}) hanno dimostrato che
ogni 4-varieta' simplettica e' in modo naturale asintoticamente realizzata
tramite un tale rivestimento generico, ed hanno proposto di
utilizzare un quoziente appropriato di $\pi_1(\PP^2_{\C} - B)$ 
per produrre invarianti di strutture simplettiche.

In questo contesto (cf. \cite{a-d-k}) gli esempi di Moisezon sono stati
rivisitati in modo piu' semplice tamite operazioni di chirurgia simplettica,
mentre Auroux, Donaldson, Katzarkov e Yotov (\cite{a-d-k-y}) hanno 
verificato che gli invarianti prodotti non riescono per\'o a distinguere 
al di la' della struttura topologica.

Questo quadro emerge anche da una serie di lavori di Moisezon e Teicher
( si veda ad esempio \cite{m-t})
e l' unico suggerimento interessante, sebbene complicato, che mi sento di dare
 e' il seguente.

Consideriamo una superficie di tipo generale $S$ ed una sua 
immersione $m$-pluricanonica: una proiezione generica su $\PP^3_{\C} $ 
da' allora una superficie con una curva doppia $\Gamma '$, mentre una
proiezione successiva su $\PP^2_{\C} $ ci da' non solo la curva
di diramazione $B$, ma anche la curva $\Gamma$ immagine di $\Gamma '$.

Anche se la congettura di Chisini ci dice che da un punto di vista olomorfo
la curva $B$ determina la superficie $S$ e quindi la curva $\Gamma$,
da cio' non segue che il gruppo fondamentale $\pi_1(\PP^2_{\C} - B)$
determini il gruppo $\pi_1(\PP^2_{\C} - B - \Gamma)$.

Sarebbe molto bello poter calcolare questo secondo gruppo
fondamentale, annche se solo in casi particolari.

{\bf Acknowledgements.}
Vorrei ringraziare sentitamente Fabio Tonoli per utili consigli
e per avere prodotto i file eps delle figure.

\noindent
{\bf Author's address:}

\bigskip

\noindent 
Prof. Fabrizio Catanese\\
Lehrstuhl Mathematik VIII\\
Universit\"at Bayreuth, NWII\\
 D-95440 Bayreuth, Germany

e-mail: Fabrizio.Catanese@uni-bayreuth.de

\end{document}